\documentclass[11pt]{amsart}

\usepackage{amsmath,amsfonts,amsthm,amssymb, mathtools}
\usepackage{setspace}
\usepackage{ mathrsfs }
\usepackage{ marvosym }
\usepackage{stackengine}
\usepackage{graphicx}
\usepackage{bm}
\usepackage{caption}
\usepackage{subfig}
\usepackage{floatrow}
\usepackage{tikz-cd}
\usepackage{cite}

\usepackage{Tabbing}
\usepackage{fancyhdr}
\usepackage{lastpage}
\usepackage{extramarks}

\usepackage{chngpage}
\usepackage{soul,color}
\usepackage{graphicx,float,wrapfig}
\usepackage{fullpage}
\usepackage{hyperref}
\usepackage[all]{xy}

\usepackage{enumerate}
\usepackage{marvosym}
\usepackage{enumitem}
\newcommand{\be}{\begin{enumerate}}
\newcommand{\ee}{\end{enumerate}}
\newcommand{\beq}{\begin{equation}}
\newcommand{\eeq}{\end{equation}}
\usepackage{verbatim}

\newcommand{\ba}{\begin{align*}}

\newcommand{\ea}{\end{align*}}

\newcommand{\D}{{\mathbb D}}

\newcommand{\T}{{\mathcal{T}}}

\newcommand{\C}{{\mathbb C}}

\newcommand{\ol}[1]{\overline{#1}}

\newcommand{\ra}{\rightarrow}
\newcommand{\tra}\twoheadrightarrow
\newcommand{\hra}{\hookrightarrow}

\newcommand{\tla}\twoheadleftarrow

\newcommand{\lp}{\left(}

\newcommand{\rp}{\right)}
\newcommand{\lpi}{\left|}
\newcommand{\rpi}{\right|}

\def\bea#1\eea{\begin{align*}#1\end{align*}}
\def\bc#1\ec{}

\newtheorem{Theorem}{Theorem}[section]

\newtheorem{Lemma}[Theorem]{Lemma}
\newtheorem{Proposition}[Theorem]{Proposition}
\newtheorem{Corollary}[Theorem]{Corollary}
\newtheorem{Conjecture}[Theorem]{Conjecture}

\numberwithin{equation}{section}

\let\lim\relax \DeclareMathOperator*\lim{lim\vphantom{p}}

\title{Isometric submersions of Teichm\"uller spaces are forgetful}
\date{}
\author{Dmitri Gekhtman and Mark Greenfield}

\begin{document}
\captionsetup[subfigure]{labelfont=rm} 
\captionsetup[subfloat]{captionskip = 5pt}
\captionsetup{font=footnotesize}

\maketitle

\begin{abstract}
We study the class of holomorphic and isometric submersions between finite-type Teichm\"uller spaces. 
We prove that, with potential exceptions coming from low-genus phenomena, any such map is a forgetful map 
$\T_{g,n} \ra \T_{g,m}$ obtained by filling in punctures. This generalizes a classical result of Royden and Earle-Kra asserting that biholomorphisms between finite-type Teichm\"uller spaces arise from mapping classes.
As a key step in the argument, we prove that any $\mathbb{C}$-linear embedding $Q(X)\hra Q(Y)$
between spaces of integrable quadratic differentials is, up to scale,
pull-back by a holomorphic map. We accomplish this step by adapting methods developed 
by Markovic to study isometries of infinite-type Teichm\"uller spaces. 
\end{abstract}

\section{Introduction}

\subsection{Biholomorphisms of Teichm\"uller spaces}
	Let $S_{g,n}$ denote the surface of genus $g$ with $n$ punctures
and let $\T_{g,n}$ denote the corresponding Teichm\"uller space.  A central theme in Teichm\"uller theory is the
interplay between the analytic structure of $\T_{g,n}$ and the topology and geometry of the underlying finite-type surface $S_{g,n}$.
This theme is exemplified by the result of Royden \cite{royden} asserting
that every biholomorphism of $\T_g$ with $g\geq 2$ arises from the action of a mapping class of $S_g$. To prove this, Royden first established that the Teichm\"uller metric is an invariant of the complex structure on $\T_g$ -- it coincides with the intrinsically defined Kobayashi metric. Thus, any biholomorphism of $\T_g$ is an isometry for the Teichm\"uller metric.
Then, by analyzing the infinitesimal properties of the Teichm\"uller norm,
Royden showed that any holomorphic isometry is induced by a mapping class.
Earle and Kra \cite{earlekra} later extended Royden's result to the finite-dimensional Teichm\"uller spaces $\T_{g,n}$. Finally, Markovic \cite{markovic} generalized to the infinite dimensional case, proving for any Teichm\"uller space of complex dimension $\geq 2$, that the biholomorphisms are induced by quasi-conformal self-maps of the underlying Riemann surface.

\subsection{Isometric submersions between finite-type Teichm\"uller spaces}

	Royden, Earle-Kra, and Markovic characterized holomorphic isometries between Teichm\"uller spaces - except in a few low-complexity cases, these are induced by identifications of the underlying surfaces. In this paper, we characterize a broader class of maps between finite-type Teichm\"uller spaces - the holomorphic and isometric submersions. Recall that a $C^1$ map between Finsler manifolds is an {\em isometric submersion} if the derivative maps the unit ball of each tangent space onto the unit ball of the target tangent space. Our {\bf main result} is that the holomorphic isometric submersions between Teichm\"uller spaces are all of geometric origin -
with some low genus exceptions, these submersions are precisely the forgetful maps $\T_{g,n} \ra \T_{g,m}$.
\begin{Theorem}[Main Result]\label{main}
Let $F:\T_{g,n} \ra \T_{k,m}$ be a holomorphic map which is an isometric submersion with respect to the 
Teichm\"uller metrics on the domain and range. Assume $(k,m)$ satisfies the following conditions: 
\begin{align}
\label{H1}  &\mkern-200mu \mbox{The type $(k,m)$ is {\em non-exceptional}: $2k+m \geq 5$.}\\	
\label{H2}  &\mkern-200mu \mbox{The genus $k$ is positive: $k\geq 1$.}
\end{align} 
Then $g=k$, $n\geq m$, and
up to pre-composition by a mapping class, 
$F:\T(S_{g,n}) \ra \T(S_{g,m})$
is the forgetful map induced by filling in the last $n-m$ punctures of 
$S_{g,n}$.    
\end{Theorem}

{\em Remark 1:}
Recall that we have isomorphisms $\T_{2,0}\cong \T_{0,6}$ and $\T_{1,2}\cong \T_{0,5}$
induced by hyperelliptic quotients.
Thus, our hypothesis on the type $({k,m})$ can be rephrased as follows -- $\T_{k,m}$ is of complex dimension at least 2 and
is not biholomorphic to a genus zero Teichm\"uller space $\T_{0,m}$.
We expect that it is possible to remove the genus condition:
\begin{Conjecture}
Any holomorphic and isometric submersion between
finite-dimensional Teichm\"uller spaces of complex dimension at least 2
is a composition of
\begin{enumerate}
\item
Forgetful maps $\T_{g,n} \ra \T_{g,m}$ with $m<n$.
\item
Mapping classes $\T_{g,n} \ra \T_{g,n}$.
\item
The isomorphisms $\T_{2,0}\cong \T_{0,6}$ and $\T_{1,2}\cong \T_{0,5}$.
\end{enumerate}
\end{Conjecture}

{\em Remark 2:}
The complex dimension 1 Teichm\"uller spaces $\T_{0,4}$, $\T_{1,0}$,
and $\T_{1,1}$ are all biholomorphic to the unit disk $\D$.
There are many isometric submersions $\T_g \ra \D$ -- the diagonal entries of
the canonical period matrix are examples (see \cite{mcmullen} Theorem 5.2).  

\subsection{Infinitesimal geometry of the cotangent space}
The proofs of Royden's theorem and its generalizations
hinge on the analysis of the infinitesimal geometry of the Teichm\"uller norm.
Fix a marked Riemann surface $X \in \T_{g,n}$.
Then the cotangent space $T^*_X\T_{g,n}$ identifies with the space $Q(X)$ of integrable holomorphic
quadratic differentials on $X$. With respect to this identification, the dual Teichm\"uller norm is the $L^1$ norm $\|\phi \| = \int_X \lpi \phi \rpi$.
Thus, a holomorphic isometry $F:\T_{g,n} \ra \T_{k,m}$ induces for each $X\in \T_{g,n}$ a bijective, $\mathbb{C}$-linear isometry of quadratic differential spaces $Q\lp F(X) \rp \ra Q(X)$. The core step in the proof of Royden's theorem is showing that, up to scale by a constant $e^{i\theta}$, any such isometry is pullback by a biholomorphism $X \ra F(X)$.

Our study of isometric submersions between Teichm\"uller spaces follows a similar tack. The key observation is that an isometric submersion induces
isometric embeddings of cotangent spaces (see Section \ref{Submersions}). 
An important step in our analysis is the following classification result, which is of independent interest.
\begin{Theorem}\label{quaddiff}
Let $X$ and $Y$ be finite-type Riemann surfaces. Let $\widehat{X}$ and $\widehat{Y}$ be the compact surfaces obtained by filling the punctures
of $X$ and $Y$. Assume the type $(k,m)$ of $X$
is non-exceptional: $2k+m \geq 5$.   
Let $T:Q(X) \hra Q(Y)$ be a $\mathbb{C}$-linear isometric embedding. 
Then there is a holomorphic map $h:\widehat{Y} \ra \widehat{X}$ and a constant $c\in \mathbb{C}$
of magnitude $\text{deg}(h)^{-1}$ so that $T = c \cdot h^*$.    
\end{Theorem}
{\em Remark:}
Suppose $X$ is of exceptional type $(k,m)$, so $2k+m \leq 4$. 
Then one of the following holds: 
\begin{enumerate}
\item
$\dim_\C Q(X)\leq 1$
\item
$(k,m)$ is $(2,0)$ or $(1,2)$,
in which case $Q(X)$ identifies naturally
with the quadratic differential space of a surface
of non-exceptional type $(0,6)$ or $(0,5)$, respectively.
\end{enumerate} 
Thus, Theorem \ref{quaddiff} amounts to a complete classification
of $\mathbb{C}$-linear isometric embeddings $Q(X)\ra Q(Y)$ for $X$ and $Y$ of finite type. 

To prove Theorem \ref{quaddiff}, we use methods developed by Markovic \cite{markovic} in his proof of the infinite-dimensional
generalization of Royden's theorem. (See also the paper of Earle-Markovic \cite{earle} and the thesis of S. Antonakoudis \cite{stergiosthesis}.)
Recall the {\em bi-canonical map} $\widehat{X} \ra \mathbb{P}Q(X)^*$ sending each 
$x\in \widehat{X}$ to the hyperplane in $Q(X)$ of quadratic differentials vanishing at $x$. 
The idea is to relate the bi-canonical images of $X$ and $Y$ using 
a result of Rudin \cite{rudin} on isometries of $L^p$ spaces.
The fact that $T:Q(X)\ra Q(Y)$ is an isometric embedding implies via Rudin's theorem that 
$T^*:\mathbb{P}Q(Y)^*\ra \mathbb{P}Q(X)^*$ carries the bi-canonical image
of $\widehat{Y}$ onto the bi-canonical image of $\widehat{X}$. So, there is a unique $h:\widehat{Y} \ra \widehat{X}$
making the following diagram commute:
$$
\begin{tikzcd}
\mathbb{P}Q(Y)^*\arrow{r}{T^*} & \mathbb{P}Q(X)^*  \\
\widehat{Y}  \arrow{u} \arrow{r}{h} & \widehat{X}  \arrow{u}
\end{tikzcd}
$$
In fact, Rudin's result
gives us more: for any $\phi\in Q(X)$, the map $h$ pushes the $\lpi T\phi \rpi$-measure on $\widehat{Y}$ to the $\lpi \phi \rpi$-measure
on $\widehat{X}$.
Thus, we obtain the following intermediate result: 
\begin{Proposition}\label{push}
Let $X$ and $Y$ be finite-type Riemann surfaces, with $X$ of non-exceptional type.
Suppose $T:Q(X)\hra Q(Y)$ is a $\mathbb{C}$-linear isometric embedding.
Then there is a holomorphic map $h:\widehat{Y} \ra \widehat{X}$ with the following property:
For any $\phi \in Q(X)$ and any measurable $K\subset \widehat{X}$, 
$$
\int_K \lpi \phi \rpi = \int_{h^{-1}(K)} \lpi T\phi \rpi .
$$ 
\end{Proposition}
We then use Proposition \ref{push} to derive the classification result 
Theorem \ref{quaddiff}.

\subsection{Infinitesimal to global}
The last step is to obtain the global Main Result, Theorem \ref{main},
from the infinitesimal Theorem \ref{quaddiff}. We are 
given a holomorphic and isometric submersion $F:\T_{g,n} \ra \T_{k,m}$,
with $(k,m)$ satisfying hypotheses \eqref{H1} and \eqref{H2}.
Since $(k,m)$ is assumed non-exceptional, Theorem \ref{quaddiff} gives for each
$Y \in \T_{g,n}$ a holomorphic branched cover $h_Y:\widehat{Y} \ra \widehat{F(Y)}$. 
By a dimension count, it is not the case that every Riemann surface of genus $g$
is a branched cover of a surface of genus $k$
with $1\leq k < g$. We then obtain that $g=k$.
Finally, an argument involving the universal families over
$\T_{g,n}$ and $\T_{g,m}$ shows that the map $h_Y:Y \ra F(Y)$
varies continuously in $Y\in \T_{g,n}$.
Thus, the topological type of $h_Y$ is constant in $Y$. We conclude that the map $F$ is induced by a (fixed) mapping class composed with the inclusion map on the underlying surfaces, filling in punctures. 
  
\subsection{Related work}
In this paper, we generalize Royden's theorem on isometries by studying isometric submersions between Teichm\"uller spaces.
Dually, one can attempt to generalize Royden's theorem by classifying of the holomorphic and isometric {\em embeddings}
between Teichm\"uller spaces. A claimed result of S. Antonakoudis states that the 
isometric embeddings all arise from covering constructions.

Our result on submersions complements a classic theorem of Hubbard \cite{hubbard}
asserting that there are no holomorphic sections of the forgetful map
$\T_{g,1} \ra \T_{g}$, except for the six sections in genus 2 obtained
by marking fixed points of the hyperelliptic involution.
Earle and Kra \cite{earlekra} later extended the result to the setting of forgetful maps between finite-type Teichm\"uller spaces $\T_{g,n}\ra \T_{g,m}$.
Combined, our result and the theorem of Hubbard-Earle-Kra have the following interpretation:
\begin{enumerate}
\item
Holomorphic and isometric submersions between finite-dimensional Teichm\"uller spaces are of geometric origin. (They are forgetful maps.) 
\item
These submersions do not admit holomorphic sections, unless there is a geometric reason (fixed points of elliptic involutions in genus 1
and hyperelliptic involutions in genus 2).
\end{enumerate}

We mention also a result of Antonakoudis-Aramayona-Souto \cite{aramayona} stating that any holomorphic map $\mathcal{M}_{g,n} \ra \mathcal{M}_{k,m}$
between moduli spaces is forgetful, as long as $g\geq 6$ and $k\leq 2g-2$.
One can see this is as a parallel of our result, with our metric constraint replaced by an equivariance condition. 

S. Antonakoudis \cite{stergiosthesis} was the first to study isometric submersions
in the context of Teichm\"uller theory. He proved that there is no holomorphic
and Kobayshi-isometric submersion between a finite-dimensional Teichm\"uller space and a bounded symmetric domain, provided each is of complex dimension at least two.

The classification of holomorphic isometric submersions
between bounded symmetric domains is an interesting problem. See the paper of Knese \cite{knese} for the classification of holomorphic Kobayashi-isometric submersions from the polydisk
$\mathbb{D}^n$ to the disk $\mathbb{D}$. For Teichm\"uller-theoretic applications of this class of functions on the polydisk, see \cite{gekhtman} and \cite{gmarkovic}.

\subsection{Outline}
The rest of the paper is devoted to the proofs of Theorems \ref{main}
and \ref{quaddiff}. 

Section \ref{Inf} focuses on the infinitesimal geometry of
isometric submersions between Teichm\"uller spaces. In \ref{Submersions},
we recall basic facts on isometric submersions between Finsler manifolds.
In \ref{Forget} we establish that forgetful maps between Teichm\"uller spaces
are holomorphic and isometric submersions. In \ref{Rudin}, we review a theorem of Rudin concerning isometries between $L^p$ spaces,
and in \ref{Proj}, we discuss the bi-canonical embedding $X \hra \mathbb{P}Q(X)^*$
of a Riemann surface. Then, in \ref{apply} we follow the argument of \cite{markovic} to obtain
Proposition \ref{push}. Finally, in \ref{Class}, we obtain the classification Theorem \ref{quaddiff}
of isometric embeddings between quadratic differential spaces. 

Section \ref{Glob} focuses on the global geometry of isometric submersions
$F:\T_{g,n} \ra \T_{k,m}$
and the proof of the main result, Theorem \ref{main}.
In \ref{Set}, we use Theorem \ref{quaddiff} to obtain for each $Y\in \T_{g,n}$
a non-constant holomorphic map $h_Y:\widehat{Y} \ra \widehat{F(Y)}$. 
Then we use a dimension-counting argument to show that $g=k$. In \ref{Universal}, we use properties of the universal family to show that the collection of maps 
$h_Y:Y \ra X$ varies continuously in the parameter $Y\in \T_{g,n}$.
Finally, in \ref{Complete}, we finish the proof of the main result.

\noindent {\bf Acknowledgements:} The first author would like to thank Martin M\"oller for
a helpful discussion. The second author is grateful to Lizhen Ji for raising the main questions and for helpful discussions. The second author is supported by the National Science Foundation Graduate Research Fellowship Program under Grant No. DGE\#1256260.

\section{Infinitesimal Geometry}\label{Inf}

\subsection{Isometric submersions of Finsler manifolds}\label{Submersions}
We review basic properties of isometric submersions, following \cite{finsler}.
First, we recall the relevant notion from linear algebra. An {\em isometric submersion} between normed vector spaces $V$ and $W$ is a linear map
$V \ra W$ so that the image of the closed unit ball in $V$ is the closed unit ball in $W$.
Isometric submersions and isometric embeddings of normed vector spaces are dual in
the following sense.  
\begin{Lemma}\label{Lemma}
Let $T:V \ra W$ be a linear map between normed vector spaces.
\begin{enumerate}
\item
If $T$ is an isometric submersion, then the dual map
$T^*: W^* \ra V^*$ is an isometric embedding.
\item
If $T$ is an isometric embedding, then
$T^*: W^* \ra V^*$ is an isometric submersion.
\end{enumerate}
\end{Lemma} 
The proof of the first assertion of the Lemma is elementary. 
The second assertion is a restatement of the Hahn-Banach theorem.

An \emph{isometric submersion between Finsler manifolds} $M,N$ is a $C^1$ submersion $F:M\ra N$ such that the derivative
$dF_m: T_mM \ra T_{F(m)}N$
is an isometric submersion between tangent spaces with respect to the Finsler norms, for each $m\in M$.
We will use the characterization of isometric submersions in terms of isometric embeddings of cotangent spaces.
\begin{Corollary}\label{cor}
Let $F: M \ra N$ be a $C^1$ map of Finsler manifolds. Then $F$ is an isometric submersion if and only if for each $m\in M$, the coderivative
$$dF_m^*: T^*_{F(m)}N \ra T^*_{m}M$$ is an isometric embedding of cotangent spaces with
respect to the dual Finsler norms.   
\end{Corollary}

\subsection{Forgetful maps between Teichm\"uller spaces}\label{Forget}
We recall basic properties of forgetful maps between Teichm\"uller spaces, and in particular observe that these maps are holomorphic and isometric submersions. 
Let $F:\T_{g,1} \ra \T_{g}$ be the forgetful map; for each $X \in \T_{g,1}$,~
$F(X)$ is the marked Riemann surface obtained by filling in the puncture of $X$.
The cotangent space $T^*_X\T_{g,1} = Q(X)$ consists of holomorphic quadratic differentials
on $X$ with at worst a simple pole at the puncture, while $T^*_{F(X)}\T_{g} = Q\lp  F(X) \rp$
consists of those quadratic differentials on $X$ which extend holomorphically over the puncture.
The co-derivative $dF_X^*$ is the inclusion $Q\lp F(X) \rp \hra Q\lp X \rp$,
which is clearly isometric and complex-differentiable. Thus, $F$ is a holomorphic and isometric submersion.  
The same reasoning shows that any forgetful map $\T_{g,n} \ra \T_{g,m}$
is an isometric submersion. 

\subsection{Rudin's Equimeasurability Theorem}\label{Rudin}
We will need a general result of Rudin concerning isometries between subspaces of $L^p$ spaces. 
Markovic \cite{markovic} used this result in the $p=1$ case to extend Royden's theorem to Teichm\"uller spaces of infinite dimension, 
and Earle-Markovic \cite{earle} used the result to give a new and illuminating proof of Royden's theorem in the finite-dimensional case.

\begin{Proposition}[Rudin \cite{rudin}, Theorem 1]
\label{rudinprop}
Let $p$ be a positive real number which is not an even integer. Let $X$ and $Y$ be sets with finite positive measures $\mu$ and $\nu$ respectively. Let $l$ be a positive integer. Suppose $f_1,\ldots,f_l$ in $L^p(\mu,\C)$, and $g_1,\ldots,g_l$ in $L^p(\nu,\C)$ satisfy the following condition:
\begin{equation}
\label{rudinassumption}
\int_X\bigg|1+\sum_{j=1}^l\lambda_jf_j\bigg|^pd\mu =\int_Y\bigg|1+\sum_{j=1}^l\lambda_jg_j\bigg|^pd\nu,\ \text{ for all }(\lambda_1,\ldots,\lambda_l)\in\C^l. 
\end{equation}
If $F=(f_1,\ldots,f_l)$ and $G=(g_1,\ldots,g_l)$, then the maps $F:X\to\C^l$ and $G:Y\to\C^l$ satisfy the following equimeasurability condition:
\begin{equation}
\label{equimeas}
\mu(F^{-1}(E)) = \nu(G^{-1}(E))\ \text{ for each Borel set }E\subseteq\C^l.
\end{equation}
\end{Proposition}
Equation \eqref{rudinassumption} is an assumption on the moments of the $\mathbb{C}^l$-valued random variables $F$ and $G$.
The conclusion \eqref{equimeas} is that $F$ and $G$ have the same distribution. In other words, the pushforward measures $F_*(\mu)$ and $G_*(\nu)$
on $\mathbb{C}^l$ are equal.

\subsection{Projective embeddings of Riemann surfaces}\label{Proj}
In this section, we establish the setting for our application of Rudin's theorem.
Let $L$ be a holomorphic line bundle over a compact Riemann surface $\widehat{X}$,
and let $\mathcal{O}(L)$ denote the space of holomorphic sections of $L$. 
There is a holomorphic map $\widehat{X}\ra \mathbb{P}\mathcal{O}\lp L\rp^*$
sending $x\in \widehat{X}$ to the hyperplane in $\mathcal{O}(L)$ consisting of sections which vanish at $x$.
An argument using the Riemman-Roch theorem (see \cite{Narasimhan} p. 55) shows that if the degree of $L$ is at least $2g+1$,
then the map $\widehat{X}\ra \mathbb{P}\mathcal{O}\lp L\rp^*$ is an embedding.

Now, let $X$ be a Riemann surface of type $(g,n)$. Denote by $\widehat{X}$ the compact, genus $g$ Riemann surface
obtained by filling in the punctures of $X$. The space $Q(X)$ consists of quadratic differentials which are holomorphic on $X$ and have at most
simple poles at the punctures $\widehat{X}\setminus X$. 
Thus, elements of $Q(X)$ correspond to sections of a line bundle on $\widehat{X}$ of degree $4g-4+n$.
By the preceding discussion, the associated {\em bi-canonical map} $\widehat{X} \ra \mathbb{P}Q(X)^*$ is an embedding
provided $4g-4+n \geq 2g+1$, or $2g+n\geq 5$. Thus, the surfaces $X$ of non-exceptional type are precisely those
for which $\widehat{X} \ra \mathbb{P}Q(X)^*$ is an embedding.

\subsection{Applying the equimeasurability theorem}\label{apply}
In this section, we apply the methods of \cite{markovic} to prove Proposition \ref{push}.  
We acknowledge some overlap with \cite{stergiosthesis} Section 5,
particularly in the proof of the fact that the surface $\widehat{Y}$ covers the surface $\widehat{X}$
if there is a $\mathbb{C}$-linear isometric embedding $Q(X) \hra Q(Y)$.

{\em Proof of Proposition \ref{push}: }
Let $X$ and $Y$ be Riemann surfaces of finite type. Assume $X$ is of non-exceptional type, and
denote by $\Phi:\widehat{X} \hra \mathbb{P}Q(X)^*$ the bi-canonical embedding associated to $X$. Let $T:Q(X) \ra Q(Y)$ be a $\mathbb{C}$-linear isometric embedding. 
Denote by $\Psi$ the composition $\widehat{Y} \ra \mathbb{P}Q(Y)^* \ra \mathbb{P}Q(X)^*$
of the bi-canonical map of $Y$ with the dual of $T$.
To describe the maps $\Phi$ and $\Psi$ more concretely, fix a basis $\phi_0,\ldots, \phi_k$ for $Q(X)$ and let $\psi_i = T\phi_i$
denote the images in $Q(Y)$. In terms of local coordinates $z, w$ for $\widehat{X}$ and $\widehat{Y}$, respectively, 
the maps $\Phi: \widehat{X} \ra \mathbb{P}^l$ and $\Psi: \widehat{Y} \ra \mathbb{P}^l$ are given by
$$
\Phi(z) = [\phi_0(z):\ldots: \phi_l(z)],~\Psi(w) = [\psi_0(w):\ldots: \psi_l(w)].
$$

Now, consider the rational functions $f_i = \frac{\phi_i}{\phi_0}$ on $\widehat{X}$ and $g_i = \frac{\psi_i}{\psi_0}$ on $\widehat{Y}$, with $i=1,\ldots, l$.
Form the $\mathbb{C}^l$-valued maps $F = (f_1,\ldots, f_l)$ and $G = (g_1,\ldots, g_l)$. The maps $F$ and $G$
are just $\Phi$ and $\Psi$ viewed as rational maps to $\mathbb{C}^l$.
 
Let $\mu$ denote the $|\phi_0|$-measure on $\widehat{X}$; that is, 
$$
\mu(K) = \int_K \lpi \phi_0 \rpi
$$ 
for any measurable $K\subset \widehat{X}$.
Similarly, let $\nu$ denote the $|\psi_0|$-measure on $\widehat{Y}$.
Then $f_i$ and $g_i$ are $L^1$ functions with with respect to the measures $\mu$ and $\nu$.
The assumption that $T$ is isometric and $\mathbb{C}$-linear translates precisely to the hypothesis \eqref{equimeas} of Rudin's theorem:  
 
\begin{align*}
\int_{\widehat{X}} \lpi 1 + \sum_{i=1}^l \lambda_i f_i \rpi d\mu &= \int_{\widehat{X}} \lpi \phi_0 + \sum_{i=1}^l \lambda_i\phi_i  \rpi \\
											      &= \int_{\widehat{Y}} \lpi \psi_0 + \sum_{i=1}^l \lambda_i\psi_i \rpi = \int_{\widehat{Y}} \lpi 1 + \sum_{i=1}^l \lambda_i g_i \rpi d\nu.										                                                                                                 
\end{align*}
Note that we used $\mathbb{C}$-linearity of $T$ in the second equality.
We conclude that the measures $F_*(\mu)$ and $G_*(\nu)$ on $\mathbb{C}^l$ are equal.
What amounts to the same thing,
the measures $\Phi_*(\mu)$ and $\Psi_*(\nu)$ on $\mathbb{P}^k$ are equal. 

We now show that $\Phi$ and $\Psi$ have the same image. 
To this end, note that the measure $\Psi_*(\nu) = \Phi_*(\mu)$ has as its support the compact set $\Phi(\widehat{X})$.
Since $\Psi$ is continuous and since $\nu$ assigns nonzero measure to each open set of $\widehat{Y}$, 
we conclude $\Psi(\widehat{Y}) \subset \Phi(\widehat{X})$.  
Thus, there is a unique holomorphic map $h:\widehat{Y} \ra \widehat{X}$ so that $\Psi = \Phi\circ h$.
Obviously, $\Psi$ is not constant and so neither is $h$. In particular, $h$ is a branched cover and
$\Psi(\widehat{Y})= \Phi(\widehat{X})$.

In terms of the map $h$, the equimeasurability condition $\Psi_*(\nu) = \Phi_*(\mu)$ becomes simply $h_*(\nu) = \mu$.
Thus, for any measurable $K\subset \widehat{X}$ we have
\begin{align*}
\int_K \lpi \phi_0 \rpi = \mu(K)
                            =  \nu\lp  h^{-1}(K) \rp =  \int_{h^{-1}(K) } \lpi T\phi_0 \rpi. \\
\end{align*}
Since $\phi_0$ was chosen arbitrarily, we have the desired equality
\begin{align*}
\int_K |\phi| = \int_{h^{-1}(K)} \lpi T\phi \rpi  
\end{align*}
for any $\phi \in Q(X)$ and any measurable $K\subset \widehat{X}$.
This completes the proof of Proposition \ref{push}.

\subsection{Completing the classification of isometric embeddings}\label{Class}
Let $\phi\in Q(X)$ and write $\psi = T\phi$. Proposition \ref{push} says 
\begin{equation} \label{measure}
\int_{h^{-1}(K)} \lpi \psi \rpi  = \int_K |\phi|
\end{equation}
for any measurable $K\subset \widehat{X}$.  
To complete the proof of Theorem \ref{quaddiff}, we must show that $\psi$ is a scalar multiple of the pullback $h^*\phi$. By working over an appropriate coordinate chart in $X$, we will reduce the proof to the following elementary lemma.

\begin{Lemma}\label{harmonic}
Let $g$ be a real-valued function defined on a domain in $\mathbb{C}$.
If both $g$ and $e^g$ are harmonic, then $g$ is constant.
\end{Lemma}  
{\em Proof:}
Compute 
\begin{align*}
0 = \lp e^g\rp_{z\ol{z}}
   = e^g\lp g_zg_{\ol{z}} + g_{z\ol{z}} \rp 
    =e^g g_z g_{\ol{z}}.
\end{align*}
Thus, $g$ is either holomorphic or anti-holomorphic. Since $g$ is real-valued, it follows that it is constant.
\qed

Returning to the proof of Theorem \ref{quaddiff},
fix a coordinate chart $(U,z)$ in $X$ on which $\phi = (dz)^2$.
(Recall that one achieves this by integrating a local square root of $\phi$.) 
Shrinking $U$ if necessary, assume $U$ is evenly covered by $h$
and that $\psi$ has no zeros or poles in $h^{-1}(U)$.
Write $h^{-1}(U)$ as a disjoint union of coordinate charts $(U_i,z_i)$,
with coordinate functions chosen so that $h:(U_i,z_i) \ra (U,z)$ is the identity function:
$$z(h(y)) = z_i(y)$$  
Let $\psi_i(z_i)(dz_i)^2$ denote the local expression for $\psi$ in $U_i$.
Let $K \subset U$ be measurable. Then equation \eqref{measure} yields
$$\int_K \lp \sum_{i=1}^{\deg(h)} |\psi_i(z)| \rp \lpi dz \rpi = \int_K \lpi dz \rpi.$$
Since $K$ was arbitrary, we have
$$ 
\sum_{i=1}^{\deg(h)} |\psi_i(z)| =1,
$$
identically on $U$. Recall that the absolute value of a holomorphic function of one variable is subharmonic.
So the function $$|\psi_1(z)| = 1 - \sum_{i=2}^{\deg(h)} |\psi_i(z)| $$
is both subharmonic and superharmonic. That is, $|\psi_1(z)|$ is harmonic.
But, since $\psi_1(z)$ is holomorphic and non-vanishing, $\log |\psi_1(z)|$ is
also harmonic. By Lemma \ref{harmonic}, $\psi_1(z)$ is identically equal to some constant $c$.
In other words, $$\psi = c\cdot h^*\phi$$ on the open set $U_1$ and thus on all of $X$.
Since $\phi\in Q(X)$ was arbitrary and $T:Q(X) \ra Q(Y)$ is linear, we have
$$T\phi = c\cdot h^*\phi$$
for all $\phi\in Q(X)$, with $c$ independent of $\phi$.
Since $T$ is an isometric embedding, we have $$|c| = \frac{\| \phi \|}{\| h^*\phi \|} = \deg(h)^{-1}.$$ 
This completes the proof of Theorem \ref{quaddiff}.

\section{The main result}\label{Glob}
\subsection{Set-up}\label{Set}
We begin the proof of the main result Theorem \ref{main}.
Let $F:\T_{g,n} \ra \T_{k,m}$ be a holomorphic and isometric submersion
of Teichm\"uller spaces. Assume $2k+m\geq 5$ and $k\geq 1$. 
By Corollary \ref{cor}, we have for each $Y \in \T_{g,n}$ that
the induced map of cotangent spaces $Q\lp F(Y) \rp  \ra Q(Y)$ is an isometric embedding.
Since $2k+m \geq 5$, Theorem \ref{quaddiff} tell us that the embedding is, up to scale, pull-back by a holomorphic branched cover
of compact surfaces
$$h_Y:\widehat{Y} \ra \widehat{F(Y)}.$$
We conclude in particular that every Riemann surface of genus $g$
admits a holomorphic branched cover of
a surface of genus $h$. We now use our assumption that $k\geq 1$. 
The following elementary lemma implies that $g=k$. 
\begin{Lemma}
Suppose $g\geq 2$. It is not the case
that every $X\in \T_g$ admits 
a holomorphic cover of a surface of genus $k$
with $1\leq k < g$.   
\end{Lemma}
{\em Proof:}
The proof is by a dimension comparison.
Suppose $1\leq k< g$ and
let $f:S_g \ra S_k$ be a degree $d$ branched cover.
Recall the Riemann-Hurwitz formula:
$$
2-2g = d\cdot (2-2k) - b,
$$ 
where $b$ is the total branch order of the cover.

We distinguish the cases $k=1$ and $k\geq 2$.
If $k\geq 2$, we have $\dim\T_g = 3g-3$ and $\dim\T_k = 3k-3$, so we get
$$
\dim\T_g = d\cdot \dim\T_k + \frac{3}{2}b.
$$
On the other hand, for a fixed topological type of branched cover, 
the space of surfaces in $Y\in\T_g$ which admit a holomorphic
cover $Y \ra X$ of that type has dimension at most
$$
\dim\T_k + b,
$$
which is less than $\dim \T_g$ since $g>k$ and thus $d>1$.

If $k=1$, then $\dim \T_g = \frac{3}{2}b$ and the dimension of the locus of
$X\in \T_g$ which admit a holomorphic cover of the given type is at most $b$.
Since $g > k=1$, the cover must have $b>0$ and so $b < \frac{3}{2}b = \dim \T_g$.

Thus, the locus of $X\in\T_g$ covering
a surface of genus less than $g$ and greater than 0
is a countable union of lower-dimensional subvarieties.
The lemma follows.
\qed

{\em Remark:}
The locus of $X\in \T_g$ which cover the square
torus is dense. This follows from the fact that the locus of
abelian differentials with rational period coordinates is dense in
the Hodge bundle over $\T_g$ \cite{zorich}.  

We conclude that $g=k$, so our submersion $F$
maps from $\T_{g,n}$ to $\T_{g,m}$ with $m\leq n$.
We are almost done: If $g\geq 2$, the covering
maps $h_Y: \widehat{Y} \ra \widehat{F(Y)}$ must be biholomorphisms.
If $g=1$, we know a priori only that $h_Y$ are (unbranched) holomorphic covers.
Since the pullback $h_Y^*$ sends $Q\lp F(Y) \rp$
into $Q(Y)$, each preimage of a puncture $p$ in $F(Y)$ must be a puncture
of $Y$. (Otherwise, $h_Y$ pulls a differential with a pole at $p$
back to a differential which is not in $Q(Y)$.)
Thus, $h_Y$ restricts to a map between the (potentially punctured) surfaces $Y$ and $X$.
The map $h_Y:Y \ra X$ and the markings $S_{g,n} \ra Y$, $S_{g,m} \ra X$
fit into a diagram  
$$
\begin{tikzcd}
S_{g,n} \arrow{d} \arrow{r} & S_{g,m} \arrow{d} \\
Y \arrow{r}{h_Y} & X
\end{tikzcd}.
$$

It remains to establish two facts.
\begin{enumerate}
\item
The maps $h_Y$ are biholomorphisms
in the $g=1$ case. 
\item
The isotopy class of $S_{g,n} \ra S_{g,m}$,
is independent of $Y \in \T_{g,n}$.  
\end{enumerate}
The key to establishing both is showing that the family
$h_Y:\widehat{Y} \ra \widehat{F(Y)}$ varies continuously in the variable $Y$.
To make this precise, we observe that
the maps $h_Y:\widehat{Y} \ra \widehat{F(Y)}$ fit together into a map
of universal curves $H:\mathcal{C}_{g,n} \ra \mathcal{C}_{g,m}$ covering the map $F:\T_{g,m} \ra \T_{g,m}$ of Teichm\"uller spaces:
$$
\begin{tikzcd}
\mathcal{C}_{g,n} \arrow{d} \arrow{r}{H} & \mathcal{C}_{g,m} \arrow{d} \\
\mathcal{T}_{g,n} \arrow{r}{F} & \mathcal{T}_{g,m}
\end{tikzcd}
$$
We will show in the next section that $H$ is continuous.
Recall $h_Y$ was constructed using the maps $X \ra \mathbb{P}Q(X)^*$
and $Y \ra \mathbb{P}Q(Y)^*$. We will leverage properties of the bundle of quadratic
differentials over Teichm\"uller space to prove that $H$ is in fact holomorphic.

\subsection{The universal curve and the cotangent bundle}\label{Universal}

We start by recalling the properties of the universal curve $\pi: \mathcal{C}_{g,n} \ra \T_{g,n}$.
A good reference for this material is \cite{Nag}. 

The map $\pi: \mathcal{C}_{g,n} \ra \T_{g,n}$ 
is a holomorphic submersion whose fiber over $X\in \mathcal{T}_{g,n}$ is exactly the compact Riemann surface $\widehat{X}$. The locations of the punctures are encoded by canonical holomorphic sections 
$$s_i: \mathcal{T}_{g,n} \ra \mathcal{C}_{g,n} ~ i=1,\ldots, n.$$ 
The point $s_i(X)\in \widehat{X}$ is the $i$th puncture of $X$.
Moreover, there is a canonical topological trivialization 
$$
\mathcal{F}_{g,n}:\mathcal{T}_{g,n}\times S_{g,n} \ra \mathcal{C}_{g,n}\setminus \bigcup_{i=1}^n s_i(\T_{g,n}),
$$ 
unique up to fiberwise isotopy,
so that the induced marking of each fiber
$$
S_{g,n} \ra \{X\}\times S_{g,n} \stackrel{\mathcal{F}_{g,n}}{\ra} X
$$ 
agrees with the marking defining $X$ as a point of $\T_{g,n}$.
The family $(\pi, \{s_i\}_{i=1}^n, \mathcal{F}_{g,n} )$ is universal among $n$-pointed marked holomorphic families of genus $g$ Riemann surfaces
(see \cite{Nag}).

Now, let $\mathcal{Q}_{g,n} \ra \T_{g,n}$ denote the bundle
of integrable holomorphic quadratic differentials over Teichm\"uller space.
Let $\mathbb{P}\mathcal{Q}_{g,n}^* \ra \T_{g,n}$ denote the associated holomorphic bundle
of projectivized dual spaces. The bi-canonical maps $\widehat{X} \ra \mathbb{P}Q(X)^*$ fit into a map
$$\Psi:\mathcal{C}_{g,n} \ra \mathbb{P}\mathcal{Q}_{g,n}^*$$ covering
the projections to Teichm\"uller space.
We need to show that this map of bundles is holomorphic.

\begin{Proposition}\label{holomorphic}
The fiberwise bi-canonical map $\Psi:\mathcal{C}_{g,n} \ra \mathbb{P}\mathcal{Q}_{g,n}^*$ is holomorphic.
If the type $(g,n)$ is non-exceptional, then the map is a biholomorphism onto its image.
\end{Proposition}
{\em Proof:}
Since $\pi$ is a holomorphic submersion,
$\mathcal{C}_{g,n}$ is covered by product neighborhoods
$U\times V$, with $U$ open in $\T_{g,n}$ and $V$ open in $\mathbb{C}$.
Each $U\times V$ maps biholomorphically to an open neighborhood of $\mathcal{C}_{g,n}$ by a
map commuting with the projections: 
$$
\begin{tikzcd}
U\times V \arrow{d} \arrow{r}{} & \mathcal{C}_{g,n} \arrow{d} \\
U \arrow{r}{} & \mathcal{T}_{g,n}
\end{tikzcd}
$$
Given $X\in U$, the slice $\{X\}\times V$ is a holomorphic coordinate chart for the Riemann surface $\widehat{X}$.
For this reason, the product neighborhoods $U\times V$ are called {\em relative coordinate charts}
for the family $\mathcal{C}_{g,n}$.

Recall $\mathcal{Q}_{g,n} \ra \T_{g,n}$, the bundle of integrable holomorphic quadratic
differentials over Teichm\"uller space. 
A section $q: \T_{g,n} \ra \mathcal{Q}_{g,n}$ can be thought of
as a fiberwise quadratic differential on $\mathcal{C}_{g,n}$.
In a relative coordinate chart $U\times V$, the differential $q$ takes the form  
$q(X, z)(dz)^2$. It follows by a result of Bers \cite{Bers} that a section $q:\T_{g,n} \ra \mathcal{Q}_{g,n}$
is holomorphic if and only if $(X,z)\mapsto q(X,z)$ is meromorphic in each relative
chart $U\times V$.

Now, let $U\times V$ be a relative coordinate chart for $\mathcal{C}_{g,n}$ and let $q_0,\ldots, q_k$
be a holomorphic frame for $\mathcal{Q}_{g,n} \ra \T_{g,n}$ over $U$. With respect to the choice of coordinates and frame,
the fiberwise bi-canonical map $\mathcal{C}_{g,n} \ra \mathbb{P}Q_{g,n}^*$ is expressed as the map $U\times V \ra \mathbb{P}^k$ given by
\begin{equation}\label{map}
(X,z) \mapsto [q_0(X,z):q_1(X,z):\cdots:q_k(X,z)],
\end{equation}
which is holomorphic since the $q_i(X,z)$ are meromorphic.

We conclude that $\Psi:\mathcal{C}_{g,n} \ra \mathbb{P}Q_{g,n}^*$ is holomorphic, as claimed.
If $(g,n)$ is non-exceptional, then $\Psi$ restricts to an embedding on the fibers
of $\mathcal{C}_{g,n} \ra \T_{g,n}$. Since the fibers are compact, $\Psi$ is a biholomorphism onto its 
image.
\qed

We now prove the main result of this subsection.
\begin{Proposition}
The map $H:\mathcal{C}_{g,n} \ra \mathcal{C}_{g,m}$ defined in the last section is holomorphic.
\end{Proposition}
{\em Proof:}
Consider the following diagram.
$$
\begin{tikzcd}
\mathcal{C}_{g,n} \arrow{rd} \arrow{r}{\Psi} & \mathbb{P}Q_{g,n}^* \arrow{d}\arrow{r}{F_*} & \mathbb{P}Q_{g,m}^*\arrow{d} & \arrow{l}[swap]{\Phi} \mathcal{C}_{g,m}\arrow{ld} \\
& \mathcal{T}_{g,n} \arrow{r}{F} & \mathcal{T}_{g,m} &
\end{tikzcd}
$$
Here, $\Psi$ and $\Phi$ denote the fiberwise bi-canonical maps,
which are holomorphic by Proposition \ref{holomorphic}.
The map $F_*$ can be viewed in two ways. 
\begin{enumerate}
\item
$F_*$ is the projectivization of the derivative of the holomorphic
map $F$.
\item
On the fiber over $Y\in \T_{g,n}$,
$F_*$ is the dual of the isometric embedding $dF_Y^*: Q\lp F(Y) \rp \hra Q(Y)$. 
\end{enumerate}
The first interpretation shows that $F_*$ is holomorphic.
The second interpretation, combined with the results of Section \ref{apply},
shows that $F_*\circ \Psi$ has the same image as $\Phi$. Moreover,
$H:\mathcal{C}_{g,n} \ra \mathcal{C}_{g,m}$ is the unique map so that
$$
F_*\circ \Psi = \Phi\circ H.
$$ 
But since $(g,m)$ is non-exceptional, $\Phi$ is a biholomorphism onto its image.
Thus, $H$ can be expressed as the composition of holomorphic maps
$$
\mathcal{C}_{g,n} \stackrel{F_*\circ \Psi}{\ra} \Phi\lp \mathcal{C}_{g,m}\rp
\stackrel{\Phi^{-1}}{\ra} \mathcal{C}_{g,m}.
$$
\qed

\subsection{Completing the proof of Theorem 1.1}\label{Complete}
As discussed at the end of Section \ref{Set}, each map $h_Y:\widehat{Y} \ra \widehat{X}$
sends $Y$ to $X$. Thus, there is a unique map 
$G:\T_{g,n}\times S_{g,n} \ra \T_{g,m}\times S_{g,m}$ fitting into the diagram
$$
\begin{tikzcd}
\mathcal{T}_{g,n}\times S_{g,n} \arrow{d}{\mathcal{F}_{g,n}} \arrow{r}{G} & \mathcal{T}_{g,m}\times S_{g,m} \arrow{d}{\mathcal{F}_{g,m}} \\
\mathcal{C}_{g,n} \arrow{r}{H} & \mathcal{C}_{g,m},
\end{tikzcd}
$$
where the vertical maps are the canonical trivializations discussed in the last section.
Since $H$ is continuous, the maps $S_{g,n} \ra S_{g,m}$ obtained by restricting $G$
to fibers are all isotopic. Restricting
the above commutative square to fibers, we conclude that there is a fixed $f:S_{g,n} \ra S_{g,m}$
so that 
$$
\begin{tikzcd}
S_{g,n} \arrow{d} \arrow{r}{f} & S_{g,m} \arrow{d} \\
Y \arrow{r}{h_Y} & F(Y)
\end{tikzcd}.
$$
commutes up to isotopy for all $Y \in \T_{g,n}$.
By construction, the vertical arrows are the markings
defining $Y$ and $F(Y)$ as points of Teichm\"uller space.
If $g\geq 2$, we already know that $f:S_{g,n} \ra S_{g,m}$ is one-to-one.
Thus, up to pre-composition by a mapping class, $Y\mapsto F(Y)$
is the forgetful map filling in the last $n-m$ punctures.
This completes the proof when $g\geq 2$.

To finish the proof in the case $g=1$, it suffices to establish that $f:S_{1,n} \ra S_{1,m}$
is one-to-one. We prove this by another dimension argument. The point is that, if the degree of $f$
is greater than $1$, then not every $X \in \T_{1,n}$ admits a non-constant holomorphic map to a $Y\in \T_{1,m}$.

In more detail:
Let $d$ denote the degree of the cover $S_1 \ra S_1$
obtained by extending $f$ over the punctures.
Then
$f$ factors through a degree $d$ (unbranched) cover $S_{1,dm} \ra S_{1,m}$.
$$
\begin{tikzcd}[column sep=small]
S_{1,n} \arrow{rr}{f}\arrow[dr] & & S_{1,m}  \\
& S_{1,dm} \arrow[ur]  &
\end{tikzcd}
$$ 
The covering $S_{1,dm} \ra S_{1,m}$ induces an isometric embedding
of Teichm\"uller spaces $\T_{1,m} \hra \T_{1,dm}$, while the injective map 
$S_{1,n} \ra S_{1,dm}$ induces a forgetful map $\T_{1,n} \tra \T_{1,dm}$.
These fit into the diagram
$$
\begin{tikzcd}[column sep=small]
\T_{1,n} \arrow{rr}{F} \arrow[dr, two heads] & & \T_{1,m} \arrow[dl, hook'] \\
& \T_{1,dm}  &
\end{tikzcd}
$$
Thus, $\T_{1,m} \hra \T_{1,dm}$ is surjective, which implies $d=1$.
\qed

\bibliography{submersionsbiblio}{}
\bibliographystyle{alpha}
\end{document}